\documentclass[12pt]{article}
\usepackage{amsmath,amsfonts,amssymb}
\usepackage{secdot}
\usepackage{amsthm}
\usepackage{color}
\newtheorem{theorem}{\hspace*{\parindent}Theorem}

\newcounter{theremark}

\title{Some remarks on rotation theorems
       \\  for complex polynomials
}
\author{V.N.\:Dubinin\footnote{E-mail address: \emph{dubinin@iam.dvo.ru} }
\\[10pt]
\small{\textit{Far Eastern Federal University
(FEFU), 8, Sukhanova Street, Vladivostok, 690950,
Russia}}\\\small{\textit{Institute of Applied
Mathematics, FEBRAS, 7, Radio Street, Vladivostok, 690041,
Russia}}}
\date{}

\begin{document}
\maketitle

\begin{abstract}
For any complex polynomial $P(z)=c_0+c_1z+...+c_nz^n,\;\;c_n\not=0,$ having all its zeros in the unit disk $|z|\le 1,$ we consider the behavior of the function (arg$P(e^{i\theta}))'_{\theta}$ when the real argument $\theta$  changes. We give some sharp estimates of this function involving of the values  of $P(e^{i\theta}),$ arg$P(e^{i\theta})$ or the coefficients $c_k,\;k=0,1,n-1,n.$ 

\end{abstract}

\bigskip


\bigskip


\bigskip

\section{Introduction}

Let all zeros of the polynomial \begin{equation} P(z)=c_0+c_1z+...+c_n z^n,\;c_n\not=0,
\end{equation}
lie in the unit disk $|z|\le 1,\;c_k\in {\mathbb C},\;k=0,1,...,n,\;n\ge1.$ 
 In the theory of polynomials, inequality                                                                            \begin{equation} {\rm Re}\;\frac{zP'(z)}{P(z)}\ge \frac{n}{2},
\end{equation} 
is well known, which is valid for all points $z$ on the circle $|z|= 1$ that are different from the zeros of the polynomial $P$. The equality in (2) takes place if and only if  all zeros of $P$ lie on $|z|= 1$ (see, for example, [1, p. 439]).  Inequality (2) can be interpreted as a rotation theorem for the complex polynomial  $P$ on the circle $|z|= 1$:                                                                     \begin{equation} ({\rm arg}P(z))'_{\theta}\ge \frac{n}{2},
\end{equation} 
$z=e^{i\theta},\;0\le \theta\le 2\pi.$                                   In a series of papers beginning with [2], this author developed the geometric function theory approach to inequalities for polynomials and rational functions.  This is described in detail in the survey article [6] and some subsequent publications of the author. In particular, as an example of the application of the boundary Schwarz lemma, the following strengthening of inequality (3) is given in [5]                                                                            \begin{equation} ({\rm arg}P(z))'_{\theta}\ge \frac{n}{2}+\frac{|c_n|-|c_0|}{2(|c_n|+|c_0|)}\ge \frac{n}{2},
\end{equation}                              $|z|=1,\; P(z)\not =0.$ Equality in (4) holds for polynomials $P$ with zeros lying on the unit circle $|z| = 1$ and for any $z,\;|z|=1,\; P(z)\not =0$ [5, inequality (10)].  Note that if the polynomial $P(z)$ of the form (1) has no zeros in the disk               $|z|<1$, then the zeros of the polynomial                                  $z^n P(1/z)$ lie in the disk                    $|z|\le 1$. Applying (4) to the last polynomial, we obtain the inequality
$$ ({\rm arg}P(z))'_{\theta}\le \frac{n}{2}+\frac{|c_n|-|c_0|}{2(|c_n|+|c_0|)}\le \frac{n}{2},
$$
for all $z$ on $|z|=1$ for which                    $P(z)\not =0$. Earlier [2], a weakened version of inequality (4) was established, in which the right-hand side of (4) is replaced by the quantity  
$$ \frac{n}{2}+\frac{\sqrt{|c_n|}-\sqrt{|c_0|}}{2\sqrt{|c_n|}}\ge \frac{n}{2},
$$                             
 [2, Theorem 4]. At present, the application of the boundary Schwarz lemma to inequalities for complex polynomials had raised considerable interest (see, for example, [11], [12], [14], [18] -- [23], [25] and references therein ). As regards inequality (4) directly, we note the following results. In the papers of  Govil, Kumar [10] and Rather, Dar, Iqbal [19] inequality (4) is generalized in different directions and another proofs of the above inequality are given. Gulzar, Zargar, Akhter [11], Hussain, Ahmad [12], and  Milovanovi${\acute{\rm c}}$, Mir, and Ahmad [18]  used (4) to generalize and strengthen the well-known inequalities for polynomials. In the paper of Rather, Dar, Iqbal [20] the weakened version (4) was used.
 
        In this note, we, first, refine inequality (4) taking into account the values of the polynomial $P$ on the circle $|z| = 1$.  Further, the strengthening of (4) with the involvement of the coefficients $c_1$      and          $c_{n-1}$ of the polynomial $P$ is obtained.  In addition, we give an upper bound for the left-hand side of (3) in the form of a finite increment theorem.  The proofs of all theorems are carried out in a unified manner and are based on the geometric function theory approach [6].  In the final Section 3, some remarks are made concerning the estimation of the value $({\rm arg}R(z))'_{\theta}$     for the rational function $R$.
\section{Rotation theorems}        
        
For points $z =e^{i\theta}$ different from the zeros of a polynomial $P$, we introduce the notation 
 $$\Lambda(P,z)=2({\rm arg}P(e^{i\theta}))'_{\theta}-n.$$
 In view of (3) 
$$\Lambda(P,z)\ge 0$$
on the circle $|z| = 1$. It is natural to pose the question of estimating the rotation speed $P(z)$ when $z$ rotates on the circle $|z| = 1$, depending on the value of $P (z)$ itself.  Some progress in the solution of this question is given by
\begin{theorem} If all zeros of a polynomial $P$ of the form {\rm(1)} lie in the disc $|z| \le 1$, then for any point $z$ a on the circle $|z| = 1$ that is different from zero of the polynomial $P$, inequality                                                        \begin{equation}                                                                                     \Lambda(P,z)\ge \left|(\Lambda(P,z)+1)\;\frac{\overline{c}_0P(z)}{c_nz^n\overline{P(z)}}-1\right| 
\end{equation}
holds.  Equality in 
{\rm(5)} at the point $z=1$ is attained for a polynomial $P$ of the form                                                                                   \begin{equation}                                                                                     P(z)=(z-\alpha_1)\prod_{k=2}^n(z-\alpha_k), 
\end{equation}
where $\alpha_k,\;k=1,...,n,$ are arbitrary numbers that satisfy the relations $|\alpha_1|<1,\;|\alpha_k|=1,\;\alpha_k\not=1,\;k=2,...,n.$
 \end{theorem}
{\bf Proof.} Let 
$$P(z)=c_n\prod_{k=1}^n(z-\alpha_k),$$
and let                              $\alpha_k\not=1,\;k=1,...,n.$ In [8] (see also [9, Theorem 1]), among other results, Goryainov obtained inequality                                                                                       \begin{equation}                                                                                     \left|f'(0)-\frac{1}{f'(1)}\right|\le 1-\frac{1}{f'(1)}, 
\end{equation}
which is valid for any holomorphic map of the disk $|z|<1$ into itself, normalized by conditions              $f(0)=0,\;f(1)=1$
and $f'(1)\not=\infty$.                 Here  $f(1)$      is the angular limit of $f$         when            $z\to1$, and $f'(1)$          is the angular derivative of the function        $f$ at the point $z = 1$. It is well known that $f'(1)\ge 1.$           Equality in (7) is attained for the function 
$$f^*(z)=z\frac{1-\overline{a}}{1-a}\frac{z-a}{1-\overline{a}z}$$
 for any $a,\;|a|<1.$                        Consider the Blaschke product
$$f(z)=\frac{\overline{c}_n}{c_n}\left[\prod_{k=1}^n\frac{1-\overline{\alpha}_k}{1-\alpha_k}\right]\frac{P(z)}{z^{n-1}\overline{P(1/\overline{z})}}=z\left[\prod_{k=1}^n\frac{1-\overline{\alpha}_k}{1-\alpha_k}\right]\prod_{k=1}^n\frac{z-\alpha_k}{1-\overline{\alpha}_k z}.$$ 
If $|\alpha|=1$, then 
$$\frac{z-\alpha}{1-\overline{\alpha} z}=\frac{1}{\overline{\alpha}}\cdot\frac{z-\alpha}{1/\overline{\alpha} -z}=-\frac{1}{\overline{\alpha}},$$
and for $|\alpha|<1$ the function 
$$w=\frac{z-\alpha}{1-\overline{\alpha} z}$$
 is a linear-fractional mapping of the disk $|z|<1$            onto the disk $|w|<1,$ such that the circle $| z |  = 1$ turns to the circle $| w |  = 1$.  
Thus, the function $f$ is holomorphic in the disk  $|z|<1$       and satisfies the conditions:
$$|f(z)|<1\;\mbox{when}\;|z|<1,\;f(0)=0,\;f(1)=1\;\mbox{and}\;f'(1)\ge1. $$
Taking into account the geometric meaning of the derivative, we conclude that at the points of the circle 
$| z |  = 1$, different from the zeros of the polynomial $P$, 
$$|f'(z)|=\frac{zf'(z)}{f(z)}=\frac{z^n\overline{P(1/\overline{z})}}{P(z)}\cdot\frac{z^{n-1}P'(z)\overline{P(1/\overline{z})}-P(z)[(n-1)z^{n-2}\overline{P(1/\overline{z})}-z^{n-3}\overline{P'(1/\overline{z})}]}{(z^{n-1}\overline{P(1/\overline{z})})^2}=$$
$$=\frac{zP'(z)}{P(z)}-n+1+\frac{\overline{P'(z)}}{z\overline{P(z)}}=2{\rm Re}\frac{zP'(z)}{P(z)}-n+1$$
is satisfied.
  
  Goryainov's inequality (7) applied to the function $f$ gives 
$$\left|\frac{c_0\overline{P(1)}}{\overline{c}_nP(1)}\left[2{\rm Re}\frac{P'(1)}{P(1)}-n+1\right]-1\right|\le 2{\rm Re}\frac{P'(1)}{P(1)}-n.$$
In general, if point $z_0,\;|z_0|=1,$  is different from the zeros of the polynomial $P$, then the point $z = 1$ is different from the zeros of the polynomial $\tilde{P}(z):=P(zz_0)$.  Applying what was proved above to the polynomial $\tilde{P}$, we have
$$\left|\frac{c_0\overline{P(z_0)}}{\overline{c}_n\overline{z}^n_0P(z_0)}\left[2{\rm Re}\frac{z_0P'(z_0)}{P(z_0)}-n+1\right]-1\right|\le 2{\rm Re}\frac{z_0P'(z_0)}{P(z_0)}-n.$$ 
The resulting inequality coincides with inequality (5) for               $z=z_0.$ If now $P (z)$ has the form (6), then 
$$f(z)=f^*(z),\;\mbox{where}\;a=\alpha_1.$$
Therefore, equality is attained in (7).  Hence the equality holds in (5).  This completes the proof of Theorem 1. 

Inequality (5) is stronger than (4).  Indeed, 
$$\Lambda(P,z)\ge 1-(\Lambda(P,z)+1)\left|\frac{c_0}{c_n}\right|$$
is fulfilled from (5).  Therefore,
$$\Lambda(P,z)\ge\frac{|c_n|-|c_0|}{|c_n|+|c_0|},$$
which is equivalent to inequality (4).
 
\begin{theorem}
Let all zeros of a polynomial $P$ of the form {\rm (1)} lie in the disc $| z |  \le 1$. Then for any point $z$ on the circle $| z |  = 1$, different from the zeros of the polynomial $P$, inequality 
                                                                                    \begin{equation}                                                                                     \Lambda(P,z)\ge\frac{2(|c_0|-|c_n|)^2}{|c_n|^2-|c_0|^2+|\overline{c}_n c_1-c_0\overline{c}_{n-1}|}                                                                                    \end{equation}                                                                                     
holds\footnote{In the case of $|c_0|=|c_n|,$ the expression on the right-hand side of (8) means zero.}.  Equality in {\rm (8)} is attained for polynomials $P$ with zeros on the unit circle $| z |  = 1$ for any                            $z,\;|z|=1,\;P(z)\not=0$.
\end{theorem}
{\bf Proof.} The equality $|c_0|=|c_n|$ holds if and only if all zeros of the polynomial $P$ lie on the circle $| z |  = 1$. Therefore, the case of equality in (8) is obvious.  Further, we consider that                        $|c_0|\not=|c_n|$.  We need the following inequality (7) from [4]: 
                                                                                   \begin{equation} 
|f'(z)|\ge 1+\frac{2(1-|f'(0)|)^2}{1-|f'(0)|^2+|f''(0)/2|}.
                                                                                    \end{equation}  
Here $f$ is a holomorphic map of the disk $| z |  < 1$ into itself, normalized by the conditions: $f(0)=0,\;|f'(0)|\not=1$                                                                     and $z$ is an arbitrary point on the circle $| z |  = 1$, in which the derivative $f'(z)$     exists and $|f(z)|=1$.  In [16], Mercer presented a direct proof of (9) using Rogosinski's lemma (see also [17, inequality (15)]).  Consider the Blaschke product
$$f(z)=\frac{P(z)}{z^{n-1}\overline{P(1/\overline{z})}}=\frac{c_n z}{\overline{c}_n}\prod_{k=1}^n\frac{z-\alpha_k}{1-\overline{\alpha}_k z},$$
where $P(z)=c_n\prod_{k=1}^n(z-\alpha_k).$  As above, we see that $f$ is a holomorphic map of the disk $| z |  <1$ into itself, $f(0)=0,$           and at each point $z$ on the circle $| z |  = 1$, nonzero of the polynomial $P$, $|f(z)|=1,$ and there is a derivative  $f'(z)$.  Moreover, 
$$f'(0)=\frac{c_n}{\overline{c}_n}\prod_{k=1}^n(-\alpha_k)=\frac{c_0}{\overline{c}_n},$$
so $|f'(0)|\not=1.$  Note that the derivative is 
$$\left(\prod_{k=1}^n\frac{z-\alpha_k}{1-\overline{\alpha}_k z}\right)'=\left(\prod_{k=1}^n\frac{z-\alpha_k}{1-\overline{\alpha}_k z}\right)\sum_{k=1}^n\left(
\frac{1}{z-\alpha_k}+\frac{\overline{\alpha}_k}{1-\overline{\alpha}_kz}\right).
$$
In view of this and Vieta's formulas, we have 
$$f''(0)=2\frac{c_n}{\overline{c}_n}\prod_{k=1}^n(-\alpha_k)\sum_{k=1}^n\left(
\frac{1}{-\alpha_k}+\overline{\alpha}_k\right)=\frac{2c_0}{\overline{c}_n}\left[
\frac{c_1}{c_0}-\frac{\overline{c}_{n-1}}{\overline{c}_{n}}\right].$$
Substituting the found values of the derivatives into (9), taking into account the calculations of the $|f'(z)|$  in the proof of Theorem 1, we arrive at inequality (8).  The theorem is proved. 

 As noted by Merser [17, Remark 3.2] for the function   $f$ 
from (9)
$$|f''(0)|\le 2(1-|f'(0)|^2)$$ 
 is fulfilled.  This attracts
 $$|c_1\overline{c}_{n}-c_0\overline{c}_{n-1}|\le|c_{n}|^2-|c_0|^2.$$
Hence it is easy to see that inequality (8) is stronger  (4).

\begin{theorem} Let all zeros of a polynomial $P$ of degree n lie in the disc $| z |  \le 1$. Suppose that for some point $z_0$ on the circle $| z |  = 1$, the arc of this circle 
$$\gamma(z_0,\alpha):=\{z:\;|z|=1,\;|{\rm arg}z-{\rm arg}z_0|<\alpha\},\;\;0<\alpha<\pi,$$
does not contain the zeros of the polynomial $P$ and the increment of the value 
$$2\;{\rm arg}P(z)-n\;{\rm arg}z$$
along any curve on the arc $\gamma(z_0,\alpha)$  with an end point at $z_0$ does not exceed in absolute value $\beta,\;0<\beta<\pi$.  Then 
                                                                                      \begin{equation} 
\Lambda(P,z_0)\le\frac{{\rm tg}\frac{\beta}{2}}{{\rm tg}\frac{\alpha}{2}}.
                                                                                   \end{equation}  
Equality in {\rm (10)} is attained for a polynomial $P$ of the form 
                                                                                     \begin{equation}
                                                                                      P(z)=c_nz(z-\alpha_2)...(z-\alpha_n),
\end{equation}
$|\alpha_k|=1,\;\alpha_k\not=1,\;k=2,...,n,$ point $z_0=1$, any                  $\alpha,\;0<\alpha<\pi,$ for which the $\gamma(1,\alpha)$ does not contain                                                  $\alpha_k,\;k=2,...,n$  and  $\beta=\alpha$. \end{theorem} 
{\bf Proof.}  We can assume that the point $z_0 = 1$. In [7], among other results, inequality 
\begin{equation}                                                                                   
 |f'(1)|\le\frac{{\rm tg}\frac{\beta}{2}}{{\rm tg}\frac{\alpha}{2}}. \end{equation}  
was established.  Here f is a holomorphic self-mapping of the unit disk $|z| < 1$, $f(\gamma(1,\alpha))\subset\gamma(1,\beta)$                                  for some $\alpha$ and                                                     $\beta$, $0<\alpha<\pi,\;0<\beta<\pi,$ and $f$ has an angular limit                     $f(1)=1$ and a finite angular derivative $f'(1)$.  Equality in (12) is attained for the function $f(z)\equiv z$   and $\alpha=\beta$.                            We put as above
$$P(z)=c_n\prod_{k=1}^n(z-\alpha_k).$$ 
Let us show that the Blaschke  product
$$f(z):=\frac{\overline{c}_n}{c_n}\left[\prod_{k=1}^n\frac{1-\overline{\alpha}_k}{1-\alpha_k}\right]\frac{P(z)}{z^{n}\overline{P(1/\overline{z})}}=\left[\prod_{k=1}^n\frac{1-\overline{\alpha}_k}{1-\alpha_k}\right]\prod_{k=1}^n\frac{z-\alpha_k}{1-\overline{\alpha}_k z}.$$  
satisfies the conditions for inequality (12). Indeed, $f$ is a holomorphic function mapping the disk $| z |  <1$ into itself.  In addition, $f$ is differentiable on the arc $\gamma(1,\alpha),\;f(1)=1$                         and $|f(z)|=1$               on   $\gamma(1,\alpha)$.  In view of the hypothesis of Theorem 3, the increment of the argument of the function $f$ along any curve $\gamma$    on $\gamma(1,\alpha)$           with an end point at the point $z = 1$ does not exceed in absolute value $\beta$: 
$$|\Delta_\gamma{\rm arg} f(z)|=|\Delta_\gamma(2{\rm arg} P(z)-n\; {\rm arg}z)|\le \beta.$$
Hence,  $f(\gamma(1,\alpha))\subset\gamma(1,\beta).$                         Following the calculations carried out in the proof of Theorem 1, we see that 
$$|f'(1)|=2{\rm Re}\frac{P'(1)}{P(1)}-n.$$
Inequality (12) gives
$$2{\rm Re}\frac{P'(1)}{P(1)}-n\le\frac{{\rm tg}\frac{\beta}{2}}{{\rm tg}\frac{\alpha}{2}},$$
which is equivalent to (10) for $z_0=1$.  If now $P (z)$ has the form (11), then $f(z)\equiv z.$  In this case, $f'(1)=1$     and 
$$|\Delta_\gamma(2{\rm arg} P(z)-n\; {\rm arg}z)|=|\Delta_\gamma{\rm arg} z|\le \alpha,\;\;\gamma\subset \gamma(1,\alpha).$$
 and you can take $\beta=\alpha$.  This gives equality in (12) and hence equality in (10).  The theorem is proved.
       
\section{On rational functions}
In conclusion we consider an application of the function theory to inequalities for rational functions with prescribed poles: 
$$R(z)=\frac{P(z)}{\prod_{k=1}^n(z-a_k)},$$
where $P(z)$ is an algebraic polynomial of degree $m$ and $|a_k|>1,\;k=1,...,n.$ In some problems, extremal functions are related to the Blaschke product
$$B(z)=\prod_{k=1}^n\frac{1-\overline{a}_k z}{z-a_k},$$
which is in general defined for any system of poles $(a_1,...,a_n),\;|a_k|\not=1,\;k=1,...,n.$

The following analogue of the polynomial inequality (3) is known.  If a rational function $R$ has exactly $m$ zeros (counting multiplicities) belonging to the disk $| z |\le 1$, then       
                                                                              \begin{equation} 
({\rm arg}R(z))'_{\theta}\ge \frac{1}{2}                                                                                 (m-n+({\rm arg}B(z))'_{\theta}) 
\end{equation}                                                                                                   
for all points on the circle $| z |  = 1$ different from zeros of $R$. If all zeros of $R$ lie in the complement to the disk $| z |  <1$, then the reverse inequality is valid: 
                                                                             \begin{equation} 
({\rm arg}R(z))'_{\theta}\le \frac{1}{2}                                                                                 (m-n+({\rm arg}B(z))'_{\theta}). 
\end{equation}                                                                                                     
Equality in (13), (14) holds for $R(z)=\alpha B(z)+\beta$                         with $|\alpha|=|\beta|$.  The proof of inequalities (13) and (14) can be found in [15, Lemma 4], [3, Lemma 3 and Theorem 4], and in [3] the geometric function theory was first applied to such a range of problems. Using various versions of the boundary Schwarz lemma, Wali and Shah  [24], [26] strengthened (13), (14) in different directions. 
Kalmykov [13] recently proved two- and three-point distortion theorems for rational functions that generalize some known results on Bernstein-type inequalities for polynomials and rational functions. The rational functions under study have either majorants or restrictions on location of their zeros. The proofs in [13] are based on the new version of the Schwarz Lemma and univalence condition for  holomorphic functions. 

In this regard, it can be assumed that the application of the methods of geometric function theory will lead to new rotation theorems for rational functions with prescribed poles, in particular, to theorems generalizing Theorems 1--3 of this article.

{\bf Funding:} This work was supported by the Russian Basic Research Fund [grant number 20-01-00018].

\end{document}